\documentclass[12pt]{amsart}
\def\la{\lambda}
\def\u{\mathbf{u}}
\allowdisplaybreaks
\title[The Eigenvectors of the Right--Justified Pascal Triangle]
{The Eigenvectors of the Right--Justified Pascal Triangle:
A Shorter Proof with Generating Functions}
\author{Helmut Prodinger}

\address{ Helmut Prodinger,
Centre for Applicable Analysis and Number Theory,
 Department of Mathematics,
University of the Witwatersrand, P.~O. Wits, 
2050 Johannesburg, South Africa, email:
{\tt helmut@gauss.cam.wits.ac.za}.
}

\date{\today}

\begin{document}

\maketitle

Let $R=\Big(\binom{i-1}{n-j}\Big)_{1\le i,j\le n}$, 
$a=\frac{1+\sqrt5}{2}$,
$\la_j=(-1)^{n+j}
a^{2j-n-1}$, $1\le j\le n$, 
\begin{equation*}
u_{ij}=\sum_{k=1}^j(-1)^{i-k}\binom{i-1}{k-1}
\binom{n-i}{j-k}a^{2k-i-1}
\end{equation*}
and $\u_j=(u_{ij})_{1\le i \le n}$.

In \cite{Callan00} Callan proves that $R\u_j=\la_j\u_j$ for
$1\le j \le n$ by what he calls a {\sl bracing exercise in manipulating
binomial coefficient sums.\/} Here we give a generating function
approach that might be a bit simpler.

\bigskip

We need also the quantity $b=\frac{1-\sqrt5}{2}$.

Consider the generating function
$U_i(z)=z(1+z)^{n-i}(az+b)^{i-1}$. It is immediate that
\begin{equation*}
U_i(z)=\sum_{j=1}^nu_{ij}z^j.
\end{equation*}
We must prove that $(R\u_j)_i=(\la_j\u_j)_i$ for all $i$ with
$1\le i\le n$. Now
\begin{align*}
(R\u_j)_i&=[z^j]\sum_{k=1}^nR_{ik}U_k(z)\\
&=[z^j]\sum_{k=1}^n\binom{i-1}{n-k}z(1+z)^{n-k}(az+b)^{k-1}\\
&=[z^{j-1}](az+b)^{n-i}\sum_{k\ge0}\binom{i-1}{k}(1+z)^{k}(az+b)^{i-k-1}\\
&=[z^{j-1}](az+b)^{n-i}\big(1+b+z(1+a)\big)^{i-1}.
\end{align*}

On the other hand, 
\begin{align*}
(\la_j\u_j)_i&=[z^j]\sum_{j=1}^n(-1)^{n+j}a^{2j-n-1}u_{ij}z^j\\
&=[z^j](-1)^{n}a^{-n-1}\sum_{j=1}^nu_{ij}(-za^2)^j\\
&=[z^j](-1)^{n}a^{-n-1}U_i(-za^2)\\
&=[z^j](-1)^{n}a^{-n-1}(-za^2)(1-za^2)^{n-i}(b-a^3z)^{i-1}\\
&=[z^{j-1}](za-\tfrac1a)^{n-i}(a^2z-\tfrac ba)^{i-1},\\
\end{align*}
and the claim follows since $b=-\frac1a$, $1+a=a^2$, and $-\frac ba=1+b$.

\bibliographystyle{amsplain}


\end{document}